\documentclass[11pt]{amsart}
\title[The mapping of compact...]{The mapping of compact into the set of its Chebyshev centres is Lipschitz in the space $l^n_{\infty}$}
\author{Pyotr N. Ivanshin}
\author{Evgenii N. Sosov}
\address{Pyotr N. Ivanshin
\newline\hphantom{iii} N.G. Chebotarev RIMM 
\newline\hphantom{iii}Kazan State University
\newline\hphantom{iii}420008, Kazan, Universitetskaya, 17}

\email{Pyotr.Ivanshin@ksu.ru}

\address{Evgenii N. Sosov  
\newline\hphantom{iii} N.G. Chebotarev RIMM 
\newline\hphantom{iii} Kazan State University
\newline\hphantom{iii} 420008, Kazan, Universitetskaya, 17}

\email{Evgenii.Sosov@ksu.ru}   
\begin{document}
\maketitle
In this article the authors prove strong stability of the set of all Chebyshev centres of the bounded closed subset of the metric space. We endow the set of all compacts of the space $l^n_{\infty}$ with Hausdorff metric and prove that the map which puts in correspondence to each compact of $l^n_{\infty}$ the set of its Chebyshev centres is Lipshitz. 

\vskip10pt
\section{Notation and definitions}
\vskip10pt

We assume the notation as follows. 

$\mathbb{R}_{+}$ denotes the set of all nonnegative real numbers. 

$\mathbb{N}$ stands for the set of all natural numbers.

$(X,\rho)$ being metric space $X$ endowed with the metric $\rho$.

$B(X)$ ($B[X]$) is the set of all nonempty bounded (bounded and closed) subsets of the metric space $X$. 

$K(X)$ stands for the set of all compact subsets of $(X,\rho)$.

$|xy| = \rho (x,y)$, $|xZ| = \inf \{|xu| : u \in Z \}$, $Mx = \sup \{|ux| : u \in M\}$
for $x, \,y \in X$, $Z \subset X$, $M \in B(X)$.

$(l^n_{\infty},\|\cdot\|_{\infty})$ --- Banach space over the field of real numbers endowed with the norm

$\|<x_1,\ldots,x_n>\|_{\infty} = \max \{|x_1|,\ldots,|x_n|\}$.

$\mathbb{Pd}(\mathbb{R}^{n})$ denotes the set of all parallelepipeds that lie in $k$--planes ($1 \leq k \leq n -1$) of the space $\mathbb{R}^{n}$, the edges of  the parallelepipeds are assumed to be parallel to coordinate axes. 

$B[x,r]$, $B(x,r)$ and $(S(x,r))$ denote respectively  closed, open ball and sphere centered in the point 
$x \in (X,\rho)$ of the radius $r \geq 0$.
 
$\Sigma_{N}$ --- set of all nonempty sub sets of $(X,\rho)$, 
consisting of no more than $N$ points. Element of the set
$\Sigma_{N}$ is called $N$--net [1].

$\alpha : B[X] \times  B[X] \rightarrow \mathbb R_{+}$, 
$\alpha  (M,T) = \max \{\sup \{|xT| :  x \in M\},\sup \{|tM| : t  \in T\}\}$ ---
Hausdorff metric on the set $B[X]$ (pseudometric on the set $B(X)$) ([2], pg. 223). 

$S(N)$ is the permutation group of $N$ elements.

$X^{N} /\!\!\sim$ denotes factor-space of the space 
$X^{N}$ with respect to the equivalence relation:
$(x_1,\ldots ,x_{N}) \sim (y_1,\ldots ,y_{N})$ in case there exists $\sigma \in S(N)$ such that

$y_1=x_{\sigma (1)},\ldots ,y_{N}=x_{\sigma (N)}$. 

$\hat \alpha  : X^{N} /\!\!\sim \times  X^{N} /\!\!\sim \, \rightarrow R_{+}$, 

$\hat \alpha : ([(x_1,\ldots ,x_{N})],[(y_1,\ldots ,y_{N})])=
\min \{\max [|x_1y_{\sigma (1)}|,\ldots ,|x_{N}y_{\sigma (N)}|] : \sigma \in S(N)\}$ is metric on the set $X^{N} /\!\!\sim$ [3]. Now using the bijection $f : X^{N} /\!\!\sim \, \rightarrow \Sigma_N$, $f([(x_1,\ldots ,x_{N})]) =
\{x_1,\ldots ,x_{N}\}$ we introduce metric $\hat \alpha$ also on the space  $\Sigma_N$.

$D[M]$ --- diameter of the set $M \in \Sigma_N$.

$\overline{M}$ --- closure of the set $M\subset X$. 

$\omega (x,y)$ is the middle set (assumed to be nonempty) of the interval $[x,y]$ such that 

$\omega (x,y) = \{z \in X : 2 |xz| = 2 |yz| = |xy|\}$ for $x$, $y \in X$.

In the Euclidean space it is simply the middlepoint of the interval $[x,y]$.

Let $M \in \Sigma_N$. $R(M) = \inf\{Mx : x \in X \}$ be a Chebyshev radius of the $N$-net $M$. The point $z \in X$ is called a Chebyshev center if 
$Mz = R(M)$ [1]. $\mathrm{cheb}(M)$ denots the set of all Chebyshev centers of $M$.

The mapping $f : (X,\rho) \rightarrow (X_1,\rho_1)$ is called Lipschitz one if there exists a constant $L \geq 0$ such that 
$\rho_1(f(x),f(y)) \leq L \rho (x,y)$ for all $x$, $y \in X$. Lipschitz mapping with 
Lipschitz constant 1 is called nonexpanding
([4], pg. 10).

In the norm spaces we use the following notation

$co(M)$ --- convex hull of the set $M$ (i.e. the intersection of all convex sets comprising $M$).

$[x,y]$ --- closed interval with endpoints $x, y$.

\section{Introduction}

In this article the authors investigate behavior of the set of Chebyshev centres of compacts of Euclidean space or those of the space $l^n_{\infty}$.
It is known ([5], [6]) that Chebyshev center of the nonempty bounded subset of the Euclidean or Lobachevskii space is strongly stable, i.e. the mapping 
\newline
$\mathrm{cheb} : (B(X),\alpha) \rightarrow X$, $M \mapsto \mathrm{cheb}(M)$,
here $B(X)$ denotes the set of all nonempty bounded subsets of the space $X$,
is continuous. Here we prove (theorem 2) similar property for the set of all nonempty subsets of the arbitrary metric space.
Recall now that the restriction of the mapping $\mathrm{cheb}$ to the set of all balls of the space with inner metric is nonexpanding map
[7]. It holds true also for the restriction of this mapping to the set of all $N$--nets of Euclidean line [8] as well as for the restriction of the map to the set of all $2$--nets of the space of nonpositive Busemann curvature [9]. 
At the same time the restriction of the map $\mathrm{cheb}$ to  the set of all nonempty closed convex sets of Euclidean plane is not Lipshitz even in the neighbourhood of the closed circle [10]. Note also that if dimension of Euclidean or Lobachevskii space is greater than $1$ and $N>2$ this map is not Lipschitz in the neighbourhood of the space $\Sigma_2 (X) \subset (\Sigma_N (X),\alpha)$ [8]. 
At the same time in Hilbert space the map $\mathrm{cheb} : (\Sigma_N (X)\backslash \Sigma_{N-1} (X),\alpha) \rightarrow X$, $M \mapsto \mathrm{cheb}(M)$ 
stays locally Lipschitz [8] and for any two compacts $M, W$ the inequality $|\mathrm{cheb}(M)\mathrm{cheb}(W)| \leq \sqrt{R(M) + R(W) + \alpha (M,W)} \sqrt{\alpha (M,W)}$ holds true [11]. The authors managed to prove that in
the space $X = l^n_{\infty}$ the situation simplifies as follows: the mapping $\mathrm{cheb} : (K(X),\alpha) \rightarrow (\mathbb{Pd}(X),\alpha)$ is Lipschitz with constant $2$ (theorem 2).

\section{Statements of the results}

Let us introduce precise statements of the results of the article. To do this properly we need lemma 0 
which is part of lemma 2 from [9].

{\bf Lemma 0} Let the set $\omega (p,x)$ consist of one point for any points $p, \, x, \,y$ of the metric space $(X,\rho)$ and the inequality 
$$2 |\omega (p,x)\omega (p,y)| \leq  |xy|$$ 
hold true.

Then the inequalities 
$$|\mathrm{cheb}(M)\mathrm{cheb}(Z)| \leq \alpha (M,Z) \leq |\mathrm{cheb}(M)\mathrm{cheb}(Z)| +
(D[M] + D[Z])/2$$ 
also hold true for all $M,\, Z \in \Sigma_2 (X)$.

The following definition gives us the possibility to characterise Hausdorff metric for separable metric space.

{\bf Definition} Let $(X,\rho)$ be separable metric space, $M$, $W \in B[X]$, and $M_{*}$, $W_{*}$ be sequences consisting of elements from $M$ and $W$ respectfully, such that the closure of each sequence coinsides with the respective set. Let us define a function $\tilde \alpha$ by the formula
$$
\tilde \alpha (M,W) = \inf \{\sup [|x_ny_n| : n \in \mathbb{N}] : (x_n)_{n \in \mathbb{N}} \in M_{*}, \, (y_n)_{n \in \mathbb{N}} \in W_{*}\}.
$$

It is easy to verify both correctness of this definition and metric axioms for $\tilde \alpha$. 

{\bf Theorem 1} Let $(X,\rho)$ be separable metric space. Then $\tilde \alpha (M,W) = \alpha (M,W)$ for all $M$, $W \in B[X]$.

{\bf Lemma 1} Let $(X,\rho)$ be metric space, $M$, $M_{n} \in (B(X),\alpha)$ for $n \in \mathbb{N}$,

$\alpha (M_{n},M) \rightarrow 0$ ($n\rightarrow \infty$), $\mathrm{cheb}(M) \neq \oslash$ and $\mathrm{cheb}(M_{n}) \neq \oslash$
for $n \in \mathbb{N}$. Then 
$\lim\limits_{n\to \infty}\mathrm{cheb}(M_{n}) \subset \mathrm{cheb}(M)$.

{\bf Lemma 2} Let $X$ be uniform convex Banach space and varieties $M$ and  
$W \in B(X)$ be such that $\mathrm{cheb}(M) \in \overline{co}(M)$, $\mathrm{cheb}(W) \in \overline{co}(W)$,
$B(\mathrm{cheb}(M),R(M))\cap B(\mathrm{cheb}(W),R(W)) = \oslash$. 
Then 
$$
|\mathrm{cheb}(M)\mathrm{cheb}(W)| \leq 2 \alpha (M,W).
$$
 
{\bf Theorem 2}
Let $n > 1$, $X = l^n_{\infty}$. 
Then the mapping $\mathrm{cheb} : (K(X),\alpha) \rightarrow (\mathbb{Pd}(X),\alpha)$ is also Lipschitz with constant $2$.

Now using the obvious inequality $\alpha \leq \hat \alpha$ on the set $\Sigma_N (l^n_{\infty})\times \Sigma_N (l^n_{\infty})$ we get

{\bf Corollary}. Let $n > 1$, $X = l^n_{\infty}$.
Then the mapping $\mathrm{cheb} : (\Sigma_N (X),\hat \alpha) \rightarrow (\mathbb{Pd}(X),\alpha)$ is Lipschitz with constant $2$.

{\bf Note} The trivial modification of the proofs for the two preceding statements provides us with one similar to that of theorem 2 for the space 
$X = (l_{\infty},\|\cdot\|_{\infty})$.

{\bf Proof of theorem 1} 
Let $M$, $W \in B[X]$ and assume without loss of generality that $\alpha (M,W) = \sup \{|xW| : x \in M\}$. Then
$\sup \{|xW| : x \in M\} = \inf \{\sup \{|x_nW| : n \in \mathbb{N}\} : (x_n)_{n \in \mathbb{N}} \in M_{*}\} \leq \tilde \alpha (M,W)$. Thus $\alpha \leq \tilde \alpha$. 

Let us prove the converse relation $\alpha \geq \tilde \alpha$.
Let $M^{*} = \{x_n : n \in \mathbb{N}\}$, $W^{*} = \{y_n : n \in \mathbb{N}\}$. Then for any $\varepsilon > 0$ and arbitrary $n \in \mathbb{N}$ there exist $m(n)$, $k(n) \in \mathbb{N}$ such that 
$$
|x_ny_{m(n)}| \leq \alpha (M,W) + \varepsilon, |y_nx_{k(n)}| \leq \alpha (M,W) + \varepsilon.
$$
 Let  $W^{*}_1 = \{y_{m(n)} : n \in \mathbb{N}\}$.
We construct the set $W^{*}_2$ modifying $W^{*}_1$ with the help of the following step-by-step procedure (here $n$ takes values $1,2,\ldots$):  If $y_{m(n + 1)} \in \{y_{m(1)},\ldots,y_{m(n)}\}$ then introduce $y_{m(n + 1), n+1} = y_{m(n + 1)}$ and put in correspondence to $x_{n+1}$ element $y_{m(n + 1), n+1}$ instead of $y_{m(n + 1)}$ and replace element $y_{m(n + 1)}$ by $y_{m(n + 1), n+1}$ in the set $W^{*}_1$. Thus we get a bijection between  $M^{*}$ and  $W^{*}_3$.
Let us now put in correspondence to arbitrary element $y_n \in W^{*}_3 = W^{*}\backslash W^{*}_1$ element $x_{k(n),n} = x_{k(n)}$ of the new set $M^{*}_2$. Thus we get a bijection between $M^{*}\cup M^{*}_2$ and $W^{*}_2\cup W^{*}_3$. 
It is easy to verify that for any $\varepsilon > 0$ 
$\tilde \alpha(M,W) = \tilde \alpha(M^{*}\cup M^{*}_2,W^{*}_2\cup W^{*}_3) \leq \alpha (M,W) + \varepsilon$.

Since the number $\varepsilon > 0$ can be arbitrary small we get $\tilde \alpha(M,W)  \leq \alpha (M,W)$.
This completes the proof of theorem 1.

{\bf Proof of lemma 1}

Recall first  that $|R(M_{n}) - R(M)| \leq \alpha (M_{n},M)$ for $n \in \mathbb{N}$ [12]. 

Let $(y_n)_{n \in \mathbb{N}}$ be a sequence of points $y_{n} \in \mathrm{cheb}(M_{n})$ converging to $y \in X$. Fix arbitrary $u \in M$, then there exists a sequence $(z_n)_{n \in \mathbb{N}}$, $z_{n} \in M_{n}$ such that 
$\alpha(M_{n},M) \leq |z_{n}u_{n}| \leq 2\alpha(M_{n},M)$. Now triangle inequality and definition of Chebyshev center provide us with inequality  
$|y_{n}u| \leq  |y_{n}z_{n}| + |z_{n}u| \leq R(M_{n}) + 2\alpha(M_{n},M)$.
So since $R(M_{n}) \rightarrow R(M)$, $\alpha(M_{n},M) \rightarrow 0$, $y_n \to y$ ($n\rightarrow \infty$) we get the desired inequality $|yu|\leq R(M)$.
This completes the proof of lemma 1.

{\bf Proof of lemma 2}
Note first that clearly $\mathrm{cheb}(M) \in \overline{co}(M)$ for $M \in B(X)$ in real Hilbert and 2-dimensional Banach spaces (cf. [13], [14]).
Assume without loss of the generality that $R(W) \leq R(M)$.  Then
$$
|\mathrm{cheb}(M)\mathrm{cheb}(W)| =
|\mathrm{cheb}(M)B[\mathrm{cheb}(W),R(W)]| + R(W) \leq 
$$
$$
\leq \sup \{|xB[\mathrm{cheb}(W),R(W)]| :  x \in \overline{co}(M)\}+ R(W) =
$$
$$
= \sup \{|xB[\mathrm{cheb}(W),R(W)]| :  x \in \overline{M}\}+ R(W)  \leq 2 \alpha (M,W).
$$

This completes the proof.

{\bf Proof of theorem 2}
Recall that the space $(K(X),\alpha)$ is geodesic one [7]. Let $M$, $W \in K(X)$. Then there exists $j \in \{1,\ldots,n\}$ such that 
$$
\alpha (\mathrm{cheb}(M),\mathrm{cheb}(W)) = 
\alpha (\bigcap \limits_{x \in M} B [x,R(M)],\bigcap\limits_{y \in W} B[y,R(W)]) = 
$$
$$
\alpha (pr_{j}(\bigcap \limits_{x \in M} B [x,R(M)]),pr_{j}(\bigcap\limits_{y \in W} B[y,R(W)])),
$$
here $pr_{j}$ is projection operator onto $j$-s coordinate. Note that  
$\bigcap \limits_{x \in M} B [x,R(M)] \neq \oslash$, $\bigcap\limits_{y \in W} B[y,R(W)]\neq \oslash$ and the balls of $(X,\|\cdot\|_{\infty})$ are convex; hence
$$
pr_{j}(\bigcap \limits_{x \in M} B [x,R(M)]) =
$$
$$
=\bigcap\limits_{x \in M}pr_{j}(B [x,R(M)]) = \bigcap\limits_{x \in M} [pr_{j}(x) - R(M),pr_{j}(x) + R(M)].
$$
The second expression can be analysed in similar way.
Then there exist $x$, $u \in M$, $y$, $v  \in W$ such that 

$$\bigcap\limits_{x \in M} [pr_{j}(x) - R(M),pr_{j}(x) + R(M)] =[pr_{j}(x) - R(M),pr_{j}(u) + R(M)],$$

$$\bigcap\limits_{y \in W} [pr_{j}(y) - R(W),pr_{j}(y) + R(W)] = [pr_{j}(y) - R(W),pr_{j}(v) + R(W)].$$
Now we easily obtain the following inequalities:

$\alpha (\mathrm{cheb}(M),\mathrm{cheb}(W)) =$

$\alpha ([pr_{j}(x) - R(M),pr_{j}(u) + R(M)],[pr_{j}(y) - R(W),pr_{j}(v) + R(W)] = $

$\max \{\max \{\rho((pr_{j}(x) - R(M))),[pr_{j}(y) - R(W),pr_{j}(v) + R(W)]),$

$\rho((pr_{j}(u) + R(M)),[pr_{j}(y) - R(W),pr_{j}(v) + R(W)])\},$

$\max \{\rho((pr_{j}(y) - R(W)),[pr_{j}(x) - R(M),pr_{j}(u) + R(M)]),$

$\rho((pr_{j}(v) + R(W)),[pr_{j}(x) - R(M),pr_{j}(u) + R(M)])\}\} \leq$
 
$\max \{|pr_{j}(x) - R(M) - pr_{j}(y) + R(W)|,
|pr_{j}(u) + R(M) - pr_{j}(v) - R(W)| \leq$

$\max \{|pr_{j}(x) - pr_{j}(y)|,
|pr_{j}(u) - pr_{j}(v)|\} + |R(M) - R(W)|$.

It is not hard then to achieve the inequality 
$$
\max \{|pr_{j}(x) - pr_{j}(y)|,|pr_{j}(u) - pr_{j}(v)|\} \leq \alpha (M,W).
$$
At the same time [12] implies that $|R(M) - R(W)| \leq \alpha (M,W)$.
So 

$\alpha (\mathrm{cheb}(M),\mathrm{cheb}(W)) \leq 2 \alpha (M,W)$.

The last inequality completes the proof of theorem 2.

\bigskip
\newpage
\begin{center}

{\sc Literature}
\end{center}
\medskip

1. {\bf Garkavi, A.L.}
The best possible net and the best possible cross-section of a set in a normed space //
Am. Math. Soc., Transl., II. Ser. 39, 111-132 (1964); translation from Izv. Akad. Nauk SSSR, Ser. Mat. 26, 87-106 (1962)

2. {\bf Kuratowski, K.}
Topology. Vol. I. 
New York-London: Academic Press; Warszawa: PWN-Polish Scientific Publishers. XX, (1966).

3. {\bf Fedorchuk, V.V.;  Filippov, V.V.}
Topology of hyperspaces and its applications. (Russian) //
Mat. Kibern. 1989, No.4, 48 p. (1989)

4. {\bf Burago, D.;  Burago, Yu.;  Ivanov, S.}
A course in metric geometry. 
Graduate Studies in Mathematics. 33. Providence, RI: American Mathematical Society (AMS). xiv, (2001)

5. {\bf Belobrov, P.K.}
On the Cebysev point of a system of sets //
Izv. Vyssh. Uchebn. Zaved., Mat. 1966, No.6(55), 18-24.

6. {\bf Sosov, E.N.}
The best net, the best section, and the Chebyshev center of bounded set in infinite-dimensional Lobachevskij space. // 
Russ. Math. 43, No.9, 39-43 (1999); translation from Izv. Vyssh. Uchebn. Zaved., Mat. 1999, No.9, 42-47.

7. {\bf Sosov E.~N.} On Hausdorff intrinsic metric // Lobachevskii J.
          of Math. 2001, V. 8, P. 185-189.

8. {\bf Ivanshin P. N., Sosov E. N.} Local Lipschitz property of the map which puts in correspondence to $N$--net its 
Chebyshev center 
// www.arXiv.org /math.MG /0509275

9. {\bf Sosov E.N.} On metric space of $2$--nets of the nonpositively curved space (in Russian) (O metricheskom prostranstve vseh $2$-setei prostranstva nepolojitelnoi krivisny) // Izv. Vyssh. Uchebn. Zaved., Mat. 2004, No.10,  57-60.

10. {\bf Lang U., Pavlovic B., Schroeder V.} Extensions of Lipschitz maps into
Hadamard spaces // Geom. Funct. Anal. 2000, 10, No. 6, 1527-1553.

11. {\bf Szeptycki P., van Vleck F.S.} Centers and nearest points of sets // Proc. Am. Math. Soc. 1982, 85, 27-31.

12. {\bf Sosov, E.N.}
On the best $N$-nets of bounded closed convex sets in a special metric space. //
Russ. Math. 47, No.9, 39-42 (2003); translation from Izv. Vyssh. Uchebn. Zaved., Mat. 2003, No.9, 42-45.

13. {\bf Garkavi, A.L.}
\"Uber das Cebysevsche Zentrum und die konvexe H\"ulle einer Menge // 
Usp. Mat. Nauk 19, No.6(120), 139-145.

14. {\bf Belobrov, P.K.}
Zur Frage des Chebyshev-Zentrums einer Menge (Russian)
Izv. Vyssh. Uchebn. Zaved., Mat. 1964, No.1(38), 3-9 (1964)

\end{document}